\documentclass[10pt]{amsart}
\theoremstyle{plain}
\usepackage{amsmath, amssymb, amscd, graphicx}
\newtheorem{lem}{Lemma}[section]
\newtheorem{prop}[lem]{Proposition}
\newtheorem{dfn}[lem] {Definition}

\newtheorem{tm}{Theorem}

\newcommand{\C}{{\bf{C}}}
\newcommand{\Z}{{\bf{Z}}}

\newcommand{\im}{{{\text {im} \ }}}

\input{diagrams}

\title{Note on Khovanov link cohomology}
\author{Bojan Gornik}
\address{Princeton University, Dept. of Mathematics, Princeton, NJ 08544}
\email{bgornik@math.princeton.edu}

\begin{document}

\begin{abstract}
We extend Lee's result on $sl(2)$ Khovanov cohomology of a link $L$ to
the general $sl(n)$ case: a filtered chain complex $\overline{C}(L)$
whose spectral sequence $E_2$ term equals Khovanov cohomology is exhibited. 
We also compute $\overline{C}(L)$'s cohomology: it depends only on
linking numbers of certain sublinks of $L$.
\end{abstract}

\maketitle

\section{Introduction}

  In \cite{Khovanov} Khovanov defines an invariant of a link $L\subseteq S^3$ which takes
the form of a doubly graded abelian group, and generalizes the Jones polynomial in the
sense that group's Euler characteristic with respect to one of its gradings equals
the Jones polynomial. Given a link $L$, definition of Khovanov's invariant involves defining a graded
chain complex $C(L)$ with boundary maps which preserve the grading and taking its cohomology.
In \cite{ESL2} Lee shows that working over a field and
modifying the boundary map by adding a certain part which lowers
the grading, gives cohomology which is particularly simple - for knots it is always two dimensional.
Lee's result can also be interpreted as a construction of a filtered chain complex whose spectral 
sequence has $E_2$ term equal to Khovanov cohomology and converges to the simple form described above
for knots.

  Recently, Khovanov and Rozansky \cite{Khovanov2} introduced link cohomology which
generalizes the earlier constructruction: for each $n\ge 2$ they define link cohomology $H_n(L)$
with Euler characteristic equal to polynomial $sl(n)$-quantum invariant $P_n(L)$; the
$n=2$ case is the original Jones polynomial construction.

The aim of this paper is to extend Lee's result to the general $n\ge 2$ case. The main
results are
\begin{tm}
\label{Thm:One}

  Let $L$ be a generic planar diagram of a link in $S^3$ and $n\ge 2$. There exists a
filtered chain complex $\overline{C}(L)$ whose associated spectral sequence has
the $E_2$ term equal to $H_n(L)$.

\end{tm}

and

\begin{tm}
\label{Thm:Two}

  The dimension of $\overline{C}(L)$'s (from Theorem~\ref{Thm:One}) cohomology 
$\overline{H}(L)$ equals
$n^l$ where $l$ is the number of $L$'s components. Moreover,
to each mapping
$$ \psi : \lbrace components\,\,of\,\,L\rbrace \to \Sigma_n, $$
we can assign an element ${\bf a}_\psi\in \overline{H}(L)$ which lies in the cohomological degree
$$ \sum_{ (\varepsilon_1,\varepsilon_2)\in \Sigma_n\times\Sigma_n, 
\varepsilon_1\not=\varepsilon_2} lk(\psi^{-1}(\varepsilon_1), \psi^{-1}(\varepsilon_2)). $$
All ${\bf a}_\psi$'s generate $\overline{H}(L)$.

\end{tm}

Here, $\Sigma_n$ is a set with $n$ distinct elements, to be defined later.

  The paper is organized as follows: in section 2, we describe the construction
of $\overline{C}(L)$. We rely heavily on reader's familiarity with
\cite{Khovanov2} - the construction closely mimics the Khovanov-Rozansky construction; essentially,
the only difference is the replacement of the potential $w(x)=x^{n+1}$ corresponding to
cohomology ring of ${\bf CP}^n$ with $w(x)=x^{n+1}-(n+1)\beta^n x$ corresponding to
quantum cohomology ring of ${\bf CP}^n$. In section 3, we establish the relation between
our $\overline{C}(L)$ and the Khovanov cohomology $H_n(L)$. Finally, in section 4,
we compute $\overline{C}(L)$'s cohomology $\overline{H}(L)$.

  Let's pause briefly, to establish some notation. We fix a generic planar diagram of a link
$L\subseteq S^3$ - by abuse of notation we call it $L$ as well. Also, fix $n\ge 2$ and a
nonzero scalar $\beta\in{\bf C}$. Although our constructions of various objects will depend
on $n$ and $\beta$, the notation will not explicitely reflect that. For example, what we called
$H_n(L)$ in this introduction, will be denoted simply as $H(L)$.

\section{Construction of $\overline{C}(L)$}

  In this section we construct $\overline{C}(L)$. We will be cavalier about the details - for those, 
we refer the reader to \cite{Khovanov2}.

\begin{figure}[t]
\includegraphics{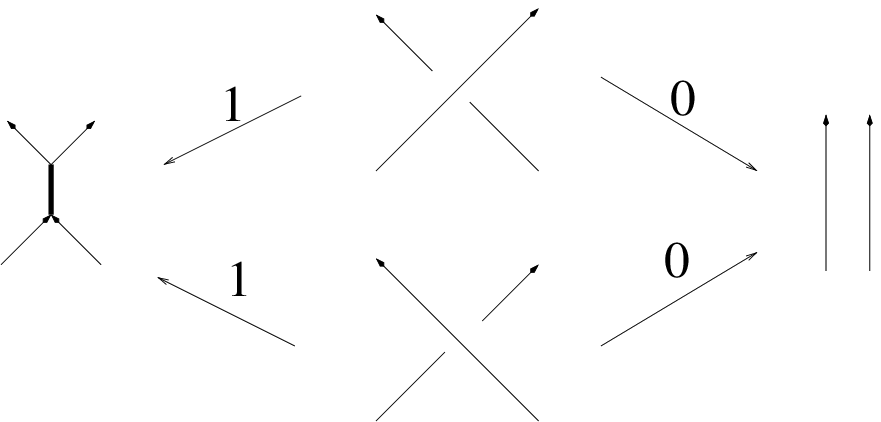}
\caption{\label{Fig:Resolutions} Resolutions of a crossing of $L$.}
\end{figure}

  Here's the outline of the construction. 
Each crossing of $L$ can be resolved in two ways, depicted in Figure~\ref{Fig:Resolutions}.
We thus obtain $2^k$ resolutions, where $k$ is the number of crossings in $L$.
To each resolution $\Gamma$ of $L$ we will assign a filtered ${\bf C}$-vector space $\overline{H}(\Gamma)$
(with some additional structure), and to a pair of resolutions $\Gamma_1$ and $\Gamma_2$,
differing locally as in Figure~\ref{Fig:Boundary}, we will assign filtered maps of degree $1$
$$\overline{\chi}_0:\overline{H}(\Gamma_0)\to \overline{H}(\Gamma_1), 
   \quad \overline{\chi}_1: \overline{H}(\Gamma_0)\to \overline{H}(\Gamma_1). $$
With these objects the filtered complex $\overline{C}(L)$ can be formed in the same manner as
in \cite{Khovanov2}, the only difference being the interpretation of the $\lbrace s\rbrace$
operator as the filtration (and not grading) shift by $s$; note that filtration 
of $\overline{H}(\Gamma)$'s induces a filtration
of $\overline{C}(L)$, and that the boundary mappings of $\overline{C}(L)$ respect this filtration, ie.
are mappings of filtered degree $0$.

  For future reference, we denote the corresponding objects in Khovanov-Rozansky construction
with $H(\Gamma), \chi_0, \chi_1, C(L)$.

  Now, we turn to the description of $\overline{H}(\Gamma)$ and maps $\overline{\chi}_0,\overline{\chi}_1$.
  
\subsection{Definition of $\overline{H}(\Gamma)$}
\label{Sbs:Hom}

  Fix a resolution $\Gamma$. We will define $\overline{H}(\Gamma)$ as the cohomology of a certain
2-periodic complex $\overline{M}(\Gamma)$. We proceed to define it.

\begin{figure}[t]
\includegraphics{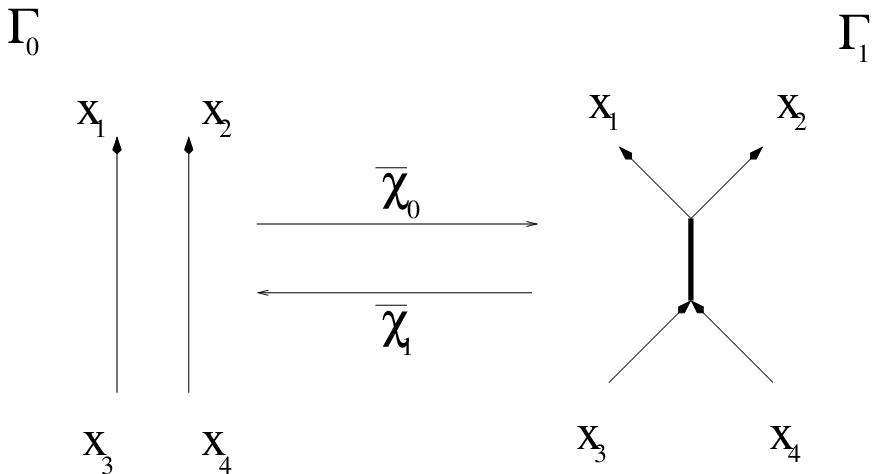}
\caption{\label{Fig:Boundary} The neighborhood of a thick edge and maps $\overline{\chi}_0,
\overline{\chi}_1$.}
\end{figure}

  First, place marks on thin edges of $\Gamma$ so that each thin edge has at least one mark
and introduce $R$, the ring of polynomials with ${\bf C}$ coefficients
in (commuting) variables $\lbrace x_i\,|\, i\in I\rbrace$ where $I$ is the set of all marks.
Introduce grading on $R$ by declaring each formal variable $x_i$ to be of degree $2$.
Then, to each arc between two neighboring marks $x, y$, oriented from $x$ to $y$, assign a factorization
$$\begin{CD} R @>{\pi_{xy} - (n+1)\beta^n}>> R\lbrace 1-n\rbrace @>{x-y}>> R, \end{CD}$$
where $\pi_{xy} = (x^{n+1}-y^{n+1})/(x-y)$.
To each thick edge (see Figure~\ref{Fig:Boundary}) assign the tensor product of factorizations
$$ \begin{CD}R\lbrace -1\rbrace @>{u_1 - (n+1)\beta^n}>> R\lbrace -n\rbrace @>{x_1+x_2-x_3-x_4}>> R\lbrace -1\rbrace 
\end{CD} $$
and
$$ \begin{CD} R @>{u_2}>> R\lbrace 3-n\rbrace @>{x_1 x_2-x_3 x_4}>> R. \end{CD} $$
Here
\begin{equation}
\label{Eqn:u1}
 u_1 = u_1(x_1,x_2,x_3,x_4) = \frac{g(x_1+x_2,x_1 x_2) - g(x_3+x_4,x_1 x_2)}{x_1+x_2-x_3-x_4}, 
\end{equation}
\begin{equation}
\label{Eqn:u2}
  u_2 = u_2(x_1,x_2,x_3,x_4) = \frac{g(x_3+x_4,x_1 x_2) - g(x_3+x_4,x_3 x_4)}{x_1 x_2 - x_3 x_4},
\end{equation}
and $g(z,w)$ is the unique two-variable polynomial for which
$$ g(x+y,xy) = x^{n+1} + y^{n+1}. $$
Finally, the (graded) 2-complex $\overline{M}(\Gamma)$ is 
defined as the tensor product (over $R$) of these factorizations.

Note that the boundary maps of $\overline{M}(\Gamma)$ are a sum of two parts: one with degree $(n+1)$ - call it
$d_0, d_1$ - and another, corresponding to the $\beta^n$ part of factorization maps,
with degree $(-n+1)$ - call it $d_0', d_1'$
$$ \begin{CD} \overline{M}^0(\Gamma) @>{d_0+d_0'}>> 
    \overline{M}^1(\Gamma) @>{d_1+d_1'}>> 
    \overline{M}^0(\Gamma).  \end{CD} $$

Therefore, the grading of $\overline{M}(\Gamma)$ does not induce a grading on its cohomology, but only a filtration.
We define the filtration $F$ of cohomology of $\overline{M}(\Gamma)$ as the one induced by the following
filtration $F'$ of $\overline{M}(\Gamma)$:
$$ F'^k\overline{M}(\Gamma) = \bigoplus_{i\le k} \overline{M}_i(\Gamma) $$
where $\overline{M}_i(\Gamma)$ is the degree $i$ direct summand of $\overline{M}(\Gamma)$.
Thus, the following is a definition of a filtered ${\bf C}$-vector space:

\begin{dfn}
Let $p(\Gamma)$ be the mod $2$ number of circles in
the modification of $\Gamma$, obtained by replacing all thick edge neighborhoods with the
$0$-resolution, ie. by performing the transformation from right to left in Figure~\ref{Fig:Boundary}
Then, 
$$ \overline{H}(\Gamma) := H^{p(\Gamma)}(\overline{M}(\Gamma)). $$
\end{dfn}

Although we will make no use of this fact, it may be worth mentioning that, while $\overline{M}(\Gamma)$'s 
${\bf Z}$ grading does not, its mod $2n$ reduction grading does descend to $\overline{H}(\Gamma)$.
Therefore, $\overline{H}(\Gamma)$ comes with a ${\bf Z}_{2n}$ grading as well. This is a consequence of the
fact that the degrees of $d_0, d_1$ are equal to those of $d_0', d_1'$ mod $2n$.

\subsection{Definition of $\overline{\chi}_0, \overline{\chi}_1$}
\label{Sbs:hi1hi0}

  The definition of maps $\overline{\chi}_0, \overline{\chi}_1$ is also analogous to the Khovanov-Rozansky construction.
Consider local reslutions $\Gamma_0, \Gamma_1$ as in Figure~\ref{Fig:Boundary}. The factorization $M_1$
associated to $\Gamma_0$ is the tensor product of factorizations, associated to the two arcs\footnote{Note that 
the corresponding formulas in \cite{Khovanov2} have indices $3$ and $4$ exchanged.}:
$$ \begin{CD}
\bmatrix R \\ R\lbrace 2-2n\rbrace \endbmatrix @>{P_0}>>
\bmatrix R \lbrace 1-n\rbrace\\ R\lbrace 1-n\rbrace \endbmatrix @>{P_1}>>
\bmatrix R \\ R\lbrace 2-2n\rbrace \endbmatrix \end{CD} $$
where
$$ P_0 = \bmatrix \pi_{13} - (n+1)\beta^n & x_2-x_4\\ \pi_{24} - (n+1)\beta^n & x_3-x_1\endbmatrix, 
\quad P_1 = \bmatrix
 x_1-x_3 & x_2-x_4\\ \pi_{24} -(n+1)\beta^n & -\pi_{13}+(n+1)\beta^n\endbmatrix, $$
$$ \pi_{ij} = (x_i^{n+1} - x_j^{n+1})/(x_i-x_j). $$

Similarly, the factorization $M_2$ associated to $\Gamma_1$ is
$$ \begin{CD}
\bmatrix R\lbrace -1\rbrace \\ R\lbrace 3-2n\rbrace \endbmatrix @>{Q_0}>>
\bmatrix R \lbrace -n\rbrace\\ R\lbrace 2-n\rbrace \endbmatrix @>{Q_1}>>
\bmatrix R\lbrace -1\rbrace \\ R\lbrace 3-2n\rbrace \endbmatrix \end{CD} $$

with
$$ Q_0 = \bmatrix u_1 - (n+1)\beta^n & x_1 x_2-x_3 x_4\\ u_2 & x_3+x_4-x_1-x_2\endbmatrix, $$
$$ Q_1 = \bmatrix
 x_1+x_2-x_3-x_4 & x_1 x_2-x_3 x_4\\ u_2 & -u_1+(n+1)\beta^n\endbmatrix. $$

  We define $\overline{\chi}_0:\overline{H}(\Gamma_0)\to\overline{H}(\Gamma_1)$ to be the map induced by the following
homomorphism of factorizations $M_1\to M_2$:
$$ U_0 : (M_1)^0 \to (M_2)^0, U_1 : (M_1)^1 \to (M_2)^1, $$
$$ U_0 = \bmatrix x_3 - x_2 & 0 \\ a_1 & 1\endbmatrix, \quad U_1 = \bmatrix x_3 & -x_2\\ -1 & 1\endbmatrix, $$
$$ a_1 = -u_2 + (u_1 + x_1 u_2 - \pi_{24})/(x_1-x_3). $$

  The map $\overline{\chi}_1:\overline{H}(\Gamma_1)\to\overline{H}(\Gamma_0)$ is induced by the following
homomorphism of factorizations $M_2\to M_1$:
$$ V_0 : (M_2)^0 \to (M_1)^0, V_1 : (M_2)^1 \to (M_1)^1, $$
$$ V_0 = \bmatrix 1 & 0 \\ -a_1 & x_3-x_2\endbmatrix, \quad V_1 = \bmatrix 1 & x_2\\ 1 & x_3\endbmatrix. $$

   With these definitions one easily checks
\begin{equation}
\label{Eqn:hi0hi1}
 \overline{\chi}_0\overline{\chi}_1 = m(x_1) - m(x_4), \quad
   \overline{\chi}_1\overline{\chi}_0 = m(x_1) - m(x_4),
\end{equation}

where $m(x_i)$ denotes the endomorphism of $\overline{H}(\Gamma_{0,1})$ induced by the multiplication by
$x_i$ endomorphism of $\overline{M}(\Gamma_{0,1})$.

\subsection{The algebra $\overline{R}(\Gamma)$; $\overline{H}(\Gamma)$ as an $\overline{R}(\Gamma)$ module}

  It will be useful for our puposes to introduce a certain (complex) algebra $\overline{R}(\Gamma)$ associated with
a resultion $\Gamma$, and the $\overline{R}(\Gamma)$ module structure of $\overline{H}(\Gamma)$. This algebra
is the analogue of $R(\Gamma)$ introduced in \cite{Khovanov3}, section 6.

\begin{dfn}
  The complex algebra $\overline{R}(\Gamma)$ is defined to be $R/({\bf a}, {\bf b})$, where 
  we view $R$, the polynomial ring in variables corresponding to all marks, as a complex algebra, 
  and ${\bf a}, {\bf b}$ are all polynomials  that appear in factorizations
  used to define $\overline{M}(\Gamma)$.
\end{dfn}

It immediatly follows, that $\overline{H}(\Gamma)$ is an $\overline{R}(\Gamma)$ module, since multiplication by
any polynomial in $({\bf a}, {\bf b})$ induces a nullhomotopic endomorphism of $\overline{M}(\Gamma)$
(see \cite{Khovanov2}, Proposition 2). Moreover, the following holds:

\begin{prop}
  As an $\overline{R}(\Gamma)$ module, $\overline{H}(\Gamma)$ is free of rank one.
\end{prop}
\begin{proof}
  This is a consequence of Proposition~\ref{Prop:Projectors} 
  and Theorems~\ref{Thm:Structure2}, \ref{Thm:Structure}.
\end{proof}

In conclusion, we give a more explicit description of $\overline{R}(\Gamma)$.

\begin{prop}
\label{Prop:rgamma}
  The algebra $\overline{R}(\Gamma)$ is spanned by generators $X_e$, where $e$ runs over all thin edges
of the resolution $\Gamma$, subject to the following relations:
for every closed circle $i$ in $\Gamma$, we have
$$ X_i^n = \beta^n, $$
and for every thick edge in $\Gamma$ the generators $X_1, X_2, X_3, X_4$ (numbered as in 
Figure~\ref{Fig:Boundary}) satisfy
$$ X_1 + X_2 - X_3 - X_4 = X_1 X_2 - X_3 X_4 = $$
$$ = u_1(X_1,X_2,X_3,X_4)-(n+1)\beta^n = u_2(X_1,X_2,X_3,X_4) = 0. $$
\end{prop}
\begin{proof}
  Since $x_i - x_j$ is in $({\bf a},{\bf b})$ for all marks $i, j$ that are adjacent and lie on the same thin edge, we
  can retain only one mark per thin edge. Also, for a closed circle with one mark $i$, the expression
  $\pi_{ii} - (n+1)\beta^n=(n+1)(x_i^n-\beta^n)$ is in $({\bf a}, {\bf b})$.
\end{proof}

Note that for any $X_e$ it holds 
$$ X_e^n = \beta^n. $$
This is evident for thin edges $e$ that correspond to closed circles. For other thin edges this follows from
thick edge relations:
$$ (n+1)(X_i^n-\beta^n) = \partial_i u_1 (X_1+X_2-X_3-X_4) + (u_1-(n+1)\beta^n)\partial_i (X_1+X_2-X_3-X_4) + $$
$$ + \partial_i u_2 (X_1 X_2 - X_3 X_4) + u_2\partial_i(X_1 X_2 - X_3 X_4) = 0, $$
where $\partial_i$ is the (formal) partial derivative with respect to $X_i$ ($i=1,2,3,4$).

\subsection{The polynomial $P_n(\Gamma)$, admissible states of $\Gamma$ and $\overline{R}(\Gamma)$}
\label{Sec:States}

  The goal of this subsection is to establish a simple relation between 
$\overline{R}(\Gamma)$ and the $sl(n)$ polynomial quantum invariant $P_n(\Gamma)$.

\begin{figure}[b]
\includegraphics{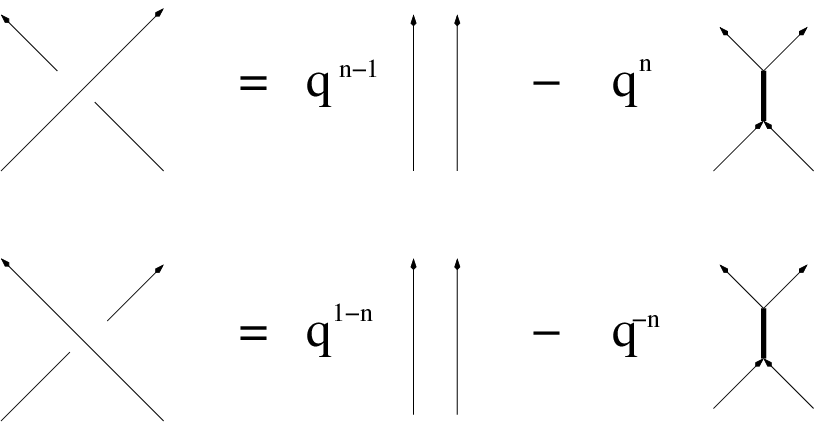}
\caption{\label{Fig:Skein} Behavior of $P_n$ with respect to resolutions.}
\end{figure}

  First, let's digress briefly to define $P_n(\Gamma) \in \Z\lbrack q,q^{-1}\rbrack$, a polynomial
in variable $q$ with positive integer coefficients, associated to a resolution $\Gamma$. 
Extending the definition of $P_n$ to links, by the rule in Figure~\ref{Fig:Skein}, we obtain a
polynomial invariant of links; in particular, $P_2$ is the Jones polynomial, and $P_n$'s are specializations
of the HOMFLY polynomial. 

  For our purposes, the following description of $P_n(\Gamma)$ will be most useful (see \cite{Jones},
Example 1.16). By $\Sigma_n$ we denote the set of $n$ complex roots of the equation $x^n-1=0$. 
Let $e(\Gamma)$ stand for the set of all thin edges of the resolution $\Gamma$. Then, an assignment
of elements in $\Sigma_n$ to all thin edges of $\Gamma$
$$ \varphi : e(\Gamma) \to \Sigma_n $$
is called {\it a state of $\Gamma$}. A state $\varphi$ is called {\it an admissible state} if at each thick edge
the $\varphi$ assignments of neighboring thin edges are either of type 1 or 2 (see the top half of
Figure~\ref{Fig:Patterns}).
The set of all states is called $S'(\Gamma)$, and the set of all admissible states is $S(\Gamma)$. Now,
\begin{equation}
\label{Eqn:states}
 P_n(\Gamma) = \sum_{\varphi\in S(\Gamma)} q^{\alpha(\varphi)}, 
\end{equation}
where $\alpha(\varphi)$ is an integer whose precise form will not interest us - it can be found in
\cite{Jones}.

  It is our goal to find a basis indexed by admissible states of $\Gamma$
  for $\C$-vector space $\overline{R}(\Gamma)$.
  We start with the following construction of idempotents $Q_\varphi$:

\begin{prop}
\label{Prop:Projectors}
  Given a state $\varphi$, we define
 \begin{equation}
 \label{Eqn:Projectors}
  Q_\varphi := \prod_{e\in e(\Gamma)}  \frac{1}{n}\left(1+\frac{X_e}{\beta\varphi(e)}
   + \frac{X_e^2}{\beta^n\varphi(e)^2} + \ldots + \frac{X_e^{n-1}}{\beta^{n-1}\varphi(e)^{n-1}}\right) 
   \in \overline{R}(\Gamma).
 \end{equation}
  It holds 
  $$ Q_\varphi^2 = Q_\varphi, \quad\forall\varphi\in S'(\Gamma), $$
  $$ Q_{\varphi_1} Q_{\varphi_2} = 0, \quad \forall\varphi_1, \varphi_2\in S'(\Gamma),\,\,\varphi_1\not=\varphi_2, $$
  $$ \sum_{\varphi\in S'(\Gamma)} Q_\varphi = 1. $$
\end{prop}
\begin{proof}
  A straightforward calculation using $X_e^n = \beta^n$.
\end{proof}

It is a matter of simple algebra to check the validity of
\begin{lem}
\label{Lem:thelemma}
  Let $\lambda_1,\lambda_2,\lambda_3,\lambda_4\in\Sigma_n$. Then 
  $$ \lambda_1+\lambda_2-\lambda_3-\lambda_4 = 0, \quad \lambda_1\lambda_2 - \lambda_3\lambda_4=0, $$
  $$  u_1(\beta\lambda_1,\beta\lambda_2,\beta\lambda_3,\beta\lambda_4) = (n+1)\beta^n, \quad
    u_2(\beta\lambda_1,\beta\lambda_2,\beta\lambda_3,\beta\lambda_4) = 0 $$
  if and only if
\begin{equation}
\label{Eqn:admissible}
   \lbrace \lambda_1, \lambda_2 \rbrace = \lbrace \lambda_3, \lambda_4\rbrace, \quad \lambda_1\not=\lambda_2. 
\end{equation}
\end{lem}

Note that (\ref{Eqn:admissible}) is precisely the condition of admissibility of a state from the top half of
Figure~\ref{Fig:Patterns}.
We now state the main result of this subsection, a structure theorem for $\overline{R}(\Gamma)$:

\begin{tm}
\label{Thm:Structure2}
   For non-admissible states $\varphi\in S'(\Gamma)-S(\Gamma)$ it holds
     $$ Q_\varphi = 0. $$
   For admissible states $\varphi\in S(\Gamma)$, however, we have
     $$  0\not=\C Q_\varphi = \overline{R}(\Gamma) Q_\varphi \quad\Longrightarrow 
         \quad\dim_{\C} \overline{R}(\Gamma)Q_\varphi = 1. $$
 Therefore, we have a direct sum decomposition of $\C$-algebras
   $$ \overline{R}(\Gamma) = \bigoplus_{\varphi\in S(\Gamma)} \C Q_\varphi. $$
 In particular,
   $$ \dim_{\C} \overline{R}(\Gamma) = P_n(\Gamma)\vert_{q=1}. $$
\end{tm}
\begin{proof}
  To prove the first statement, pick $\varphi\in S'(\Gamma)-S(\Gamma)$ and choose a thick edge where the admissibility
  condition (Figure~\ref{Fig:Patterns}) is violated. Call the neighboring thin edges $e_1,e_2,e_3,e_4$ and
  the variables corresponding to these thin edges $X_1, X_2, X_3, X_4$. 
  Since $X_e Q_\varphi = \beta\varphi(e) Q_\varphi$ for all thin edges of $\Gamma$, we have
$$   0=(X_1+X_2-X_3-X_4) Q_\varphi = \beta(\varphi(e_1)+\varphi(e_2)-\varphi(e_3)-\varphi(e_4))Q_\varphi, $$
$$   0=(X_1 X_2-X_3 X_4) Q_\varphi = \beta^2(\varphi(e_1) \varphi(e_2)-\varphi(e_3) \varphi(e_4))Q_\varphi, $$
$$   0=(u_1(X_1,X_2,X_3,X_4)-(n+1)\beta^n) Q_\varphi  = $$
$$  =  (u_1(\beta\varphi(e_1),\beta\varphi(e_2),\beta\varphi(e_3),\beta\varphi(e_4)) - (n+1)\beta^n)Q_\varphi, $$
$$   0=u_2(X_1,X_2,X_3,X_4) Q_\varphi = 
    u_2(\beta\varphi(e_1),\beta\varphi(e_2),\beta\varphi(e_3),\beta\varphi(e_4)) Q_\varphi. $$

  According to Lemma~\ref{Lem:thelemma}, one of the rightmost expressions is a nontrivial multiple of $Q_\varphi$,
  implying $Q_\varphi=0$.
  
To prove the second statement, we first observe that $Q_\varphi\not=0$. For, if $Q_\varphi$ were $0$, its 
representative in $C[X_e\,|\,e\in e(\Gamma)]$ would be in the ideal generated by all
relations from Proposition~\ref{Prop:rgamma}. But all these relations, according to
Lemma~\ref{Lem:thelemma}, give $0$ when evaluated at $X_e=\beta\varphi(e)$, whereas the expression in
(\ref{Eqn:Projectors}) clearly evaluates to $1$ at $X_e=\beta\varphi(e)$.

Introduce $P(\Gamma) = \C[X_e\,|\,e\in e(\Gamma)]/(X_e^n=\beta^n\,|\,e\in e(\Gamma))$ and elements
$Q_\varphi' \in P(\Gamma)$ given by the same expression as in (\ref{Eqn:Projectors}). Clearly, proposition
~\ref{Prop:Projectors} holds for elements $Q_\varphi'$ as well and they form a basis for
$P(\Gamma)$. Therefore, $\C Q_\varphi' = P(\Gamma) Q_\varphi'$ and the same relation holds in $\overline{R}(\Gamma)$
which is a quotient of $P(\Gamma)$.

  The last equality is a consequence of (\ref{Eqn:states}).
\end{proof}

\section{Relation to Khovanov-Rozansky construction}

  First, we briefly recall Khovanov-Rozansky's construction
of the complex $C(L)$ and state some of its properties that are relevant for us. Then, we give a proof of
Theorem~\ref{Thm:One}.

\subsection{Khovanov-Rozansky construction}

  Simply put, objects $M(\Gamma), H(\Gamma), \chi_{0,1}, C(L)$ are defined as in
subsections~\ref{Sbs:Hom}, \ref{Sbs:hi1hi0} with $\beta=0$. Note that matrices used to define
$\chi_0, \overline{\chi}_0$ (and similarly for $\chi_1, \overline{\chi}_1$) don't change.
Also, the underlying graded $\C$-spaces
of 2-complexes $M(\Gamma)$ and $\overline{M}(\Gamma)$ are equal, the difference is in the boundary maps;
in particular, omitting $d_0', d_1'$ in $\overline{M}(\Gamma)$ we obtain precisely $M(\Gamma)$.
Since $d_0, d_1$ are homogeneous (of degree $(n+1)$), the $M(\Gamma)$-cohomology $H(\Gamma)$ is
a graded space. We summarize: 

\begin{prop}
\label{Prop:Tehnika}
The 2-complex $M(\Gamma)$
$$ \begin{CD} \overline{M}^0(\Gamma) @>{d_0}>> 
    \overline{M}^1(\Gamma) @>{d_1}>> 
    \overline{M}^0(\Gamma) 
   \end{CD}
$$
has nontrivial cohomology only in cohomological degree $p(\Gamma)$, the mod 2 number of circles in
the modification of $\Gamma$, obtained by replacing all thick edge neighborhoods with the
$0$-resolution, ie. by performing the transformation from right to left in Figure~\ref{Fig:Boundary}
on all thick edges. This cohomology, by definition, is $H(\Gamma)$ and its ($q$-)graded dimension
is $P_n(\Gamma)$.
\end{prop}
\begin{proof}
  See \cite{Khovanov2}.
\end{proof}

\subsection{Proof of Theorem~\ref{Thm:One}}

  We start with the observation, that the filtered dimension of $\overline{H}(\Gamma)$ and the graded
dimension of $H(\Gamma)$ are equal.

\begin{prop}
\label{Prop:GradFilt}
  Let $k\in\Z$. If $H_k(\Gamma)$ stands for the degree $k$ direct summand of graded space $H(\Gamma)$, there
  exists an isomorphism of $\C$-vector spaces
  $$ \Phi_{k,\Gamma} : H_k(\Gamma) \to F^k\overline{H}(\Gamma)/F^{k-1}\overline{H}(\Gamma). $$
  In prticular, the filtered dimension of $\overline{H}(\Gamma)$ is $P_n(\Gamma)$.
\end{prop}
\begin{proof}
  Without loss of generality, we may assume $p(\Gamma)=1$. Then, the 2-complex $M(\Gamma)$ has trivial cohomology
in degree $0$, and cohomology in degree $1$ equals $H(\Gamma)$.
Now, let's define the map
  $$ \phi : (\ker d_1)_k \to F^k\overline{H}(\Gamma)/F^{k-1}\overline{H}(\Gamma) $$
where $(\ker d_1)_k$ is the degree $k$ direct summand of $(\ker d_1)$.
 Pick $\alpha\in (\ker d_1)_k$.and define
 $$ \phi(\alpha) = \alpha + \alpha_1 + \alpha_2 + \ldots, $$
where $\alpha_i\in (M^1)_{k-2ni}$ is of degree $k-2ni$.
Here's how we define $\alpha_1, \alpha_2, \ldots$. We need $(d_1+d_1')(\alpha+\alpha_1+\alpha_2+\ldots)=0$, or, 
equivalently, (consider the grading on $M^1$)
$$ d_1\alpha = 0, \quad d_1'\alpha+d_1\alpha_1 = 0, \quad d_1'\alpha_1+d_1\alpha_2 = 0, \ldots. $$
Since $d_0 d_1'\alpha = -d_0' d_1\alpha = 0$ (we used the fact that $\overline{M}$ is a $2$-complex) and $\ker d_0=\im d_1$,
we see that there exists $\alpha_1$ such that $d_1(-\alpha_1) = d_1'\alpha$.

In the next step we construct $\alpha_2$. Since $d_0 d_1'\alpha_1 = -d_0' d_1\alpha_1 = d_0' d_1'\alpha = 0$,
we again (due to $\ker d_0=\im d_1$) obtain $\alpha_2$ such that $d_1(-\alpha_2) = d_1'\alpha_1$.

  We continue with this algorithm, which eventually terminates since $M$ doesn't have arbitrarily low degrees.
  
  First, we have to check that our definition is a good one, ie. that the procedure gives an element in
$F^k\overline{H}(\Gamma)/F^{k-1}\overline{H}(\Gamma)$ that does not depend on choices of
$\alpha_1,\alpha_2,\ldots$. Suppose an alternative set of numbers $\alpha_1', \alpha_2', \ldots$ would also
do the job. Then 
$$(d_1+d_1')( (\alpha_1-\alpha_1') + (\alpha_2-\alpha_2') + \ldots) = 0$$
implies what we need, since we're quotientig out $F^{k-1}\overline{H}(\Gamma)$, and the above element lies
there - in fact it lies in $F^{k-2n}\overline{H}(\Gamma)$.

  Second, we note that $\phi$ is surjective: take an $M$-representative of an element 
$\alpha\in F^k\overline{H}(\Gamma)/F^{k-1}\overline{H}(\Gamma)$, strip away its part that lies in degrees
strictly lower than $k$, and what remains maps to $\alpha$ via $\phi$.

  It is easy to check $\ker\phi = (\im d_0)_k$ - we omit the details.
  
  Therefore $\Phi_{k,\Gamma}$, the mapping $(\ker d_1)_k/(\im d_0)_k \to F^k\overline{H}(\Gamma)/F^{k-1}
  \overline{H}(\Gamma)$ induced by $\phi$, is
an isomorphism.

\end{proof}

\noindent {\bf Remark:} We will not prove or need it, but it can be shown that cohomology of $\overline{M}(\Gamma)$
in cohomological degree $(p(\Gamma)+1)$ mod $2$ is trivial.

\begin{proof} (of Theorem~\ref{Thm:One}) We have to show that the $E_1$ term of the spectral sequence,
associated to filtered chain complex $\overline{C}(L)$, is isomorphic to the chain complex $C(L)$.
The $E_1$ term is a direct sum of
spaces $F^k\overline{H}(\Gamma)/F^{k-1}\overline{H}(\Gamma)$ where $\Gamma$ runs over all resolutions of $L$.
The boundary maps of $E_1$ are induced by $\overline{\chi}_0, \overline{\chi}_1$. In view of 
Proposition~\ref{Prop:GradFilt}, we will have shown the statement of the theorem, once we show that
maps $\chi_0,\chi_1$ (ie. the boundary maps of $C(L)$), via the identification from Proposition~\ref{Prop:GradFilt},
correspond to $\overline{\chi}_0, \overline{\chi}_1$ as maps of $F^k\overline{H}/F^{k-1}\overline{H}$:

\begin{diagram}
  H_k(\Gamma_0)  & \rTo^{\Phi_{k,\Gamma_0}} & F^k\overline{H}(\Gamma_0)/F^{k-1}\overline{H}(\Gamma_0) 
  & &   H_k(\Gamma_1)  & \rTo^{\Phi_{k,\Gamma_1}} & F^k\overline{H}(\Gamma_1)/F^{k-1}\overline{H}(\Gamma_1)  \\
  \dTo^{\chi_0}  & = & \dTo^{\overline{\chi}_0} & \quad {\rm AND} \quad
  & \dTo^{\chi_1}  & = & \dTo^{\overline{\chi}_1} \\
  H_{k+1}(\Gamma_1)  & \rTo^{\Phi_{k+1,\Gamma_1}} & F^{k+1}\overline{H}(\Gamma_1)/F^k\overline{H}(\Gamma_1)
  & & 
  H_{k+1}(\Gamma_0)  & \rTo^{\Phi_{k+1,\Gamma_0}} & F^{k+1}\overline{H}(\Gamma_0)/F^k\overline{H}(\Gamma_0)
\end{diagram}

This is a straightforward consequence of the definition of $\Phi$, and of the fact that maps $\chi_0$ and
$\overline{\chi}_0$ (and similarly for $\chi_1,\overline{\chi}_1$), as maps of the graded $\C$-space underlying
chain complexes $M, \overline{M}$, are equal and homogeneous (of degree $1$).
\end{proof}

\section{Computation of $\overline{H}(L)$}

  In this section we give a proof of Theorem~\ref{Thm:Two}. We start with a structure theorem for 
$\overline{H}(\Gamma)$, which establishes a basis for $\overline{H}(\Gamma)$ indexed by
admissible states of $\Gamma$. Recall that $\overline{H}(\Gamma)$ is an $\overline{R}(\Gamma)$ module; by abuse of
notation, the endomorphism of $\overline{H}(\Gamma)$ induced by the action of $Y\in\overline{R}(\Gamma)$
will also be denoted by $Y$; also, we adopt the notation introduced in subsection~\ref{Sec:States}.

\begin{figure}[t]
\includegraphics{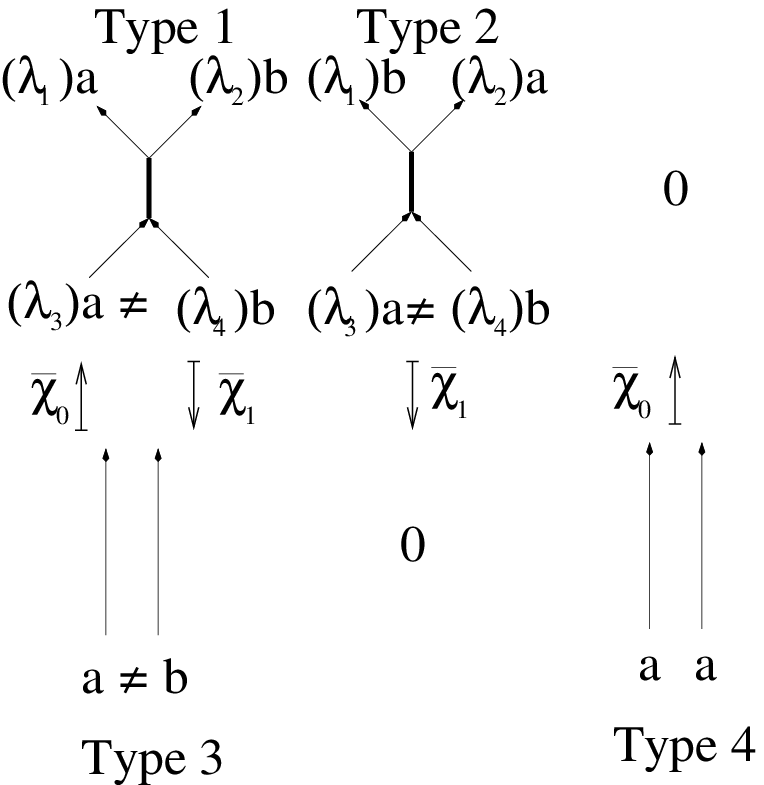}
\caption{\label{Fig:Patterns} Admissible states condition (top half). 
Action of maps $\overline{\chi}_0, \overline{\chi}_1$.}
\end{figure}

\begin{tm}
\label{Thm:Structure}
For each admissible state $\varphi\in S(\Gamma)$ the $\C$-space $Q_\varphi\overline{H}(\Gamma)$ is 
one-dimensional and
\begin{equation}
\label{Eqn:DSumDecomp}
   \overline{H}(\Gamma) = \bigoplus_{\varphi\in S(\Gamma)} Q_\varphi\overline{H}(\Gamma). 
\end{equation}
For all $x\in \overline{H}(\Gamma)$ we have
\begin{equation}
\label{Eqn:TheOne}
 x\in Q_\varphi\overline{H}(\Gamma) \Longleftrightarrow    X_e x = \beta \varphi(e) x, \quad \forall e\in e(\Gamma).
\end{equation}
\end{tm}
\begin{proof}
   The implication $\Rightarrow$ of (\ref{Eqn:TheOne}) is an obvious consequence of the definition of $Q_\varphi$;
 the implication $\Leftarrow$ follows from $x=Q_\varphi x$.
 Also, the direct sum decomposition (\ref{Eqn:DSumDecomp}) is a consequence of Proposition~\ref{Prop:Projectors} and
 $Q_\varphi=0$ for non-admissible states (Theorem~\ref{Thm:Structure2}).
   
  To finish the proof, we have to show the spaces $Q_\varphi\overline{H}(\Gamma)$ for admissible states
$\varphi$ are one-dimensional. Since
the number of admissible states of $\Gamma$ is $P_n(\Gamma)\vert_{q=1}$, and the complex dimension of $\overline{H}(\Gamma)$
equals $P_n(\Gamma)\vert_{q=1}$ (Proposition~\ref{Prop:GradFilt}), the decomposition (\ref{Eqn:DSumDecomp})
implies that it suffices to show 
\begin{equation}
  \label{Eqn:WhatWeShow}
   \dim_{\C} Q_\varphi\overline{H}(\Gamma)\ge 1. 
\end{equation}
To that end, given an admissible state $\varphi$, we will construct a nonzero element of $Q_\varphi\overline{H}(\Gamma)$.

\begin{figure}[t]
\includegraphics{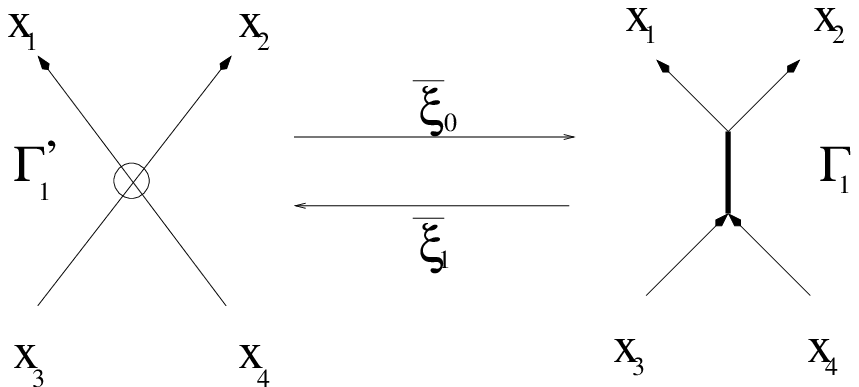}
\caption{\label{Fig:Boundary2} The maps $\overline{\xi}_0, \overline{\xi}_1$.}
\end{figure}

\begin{figure}[b]
\includegraphics{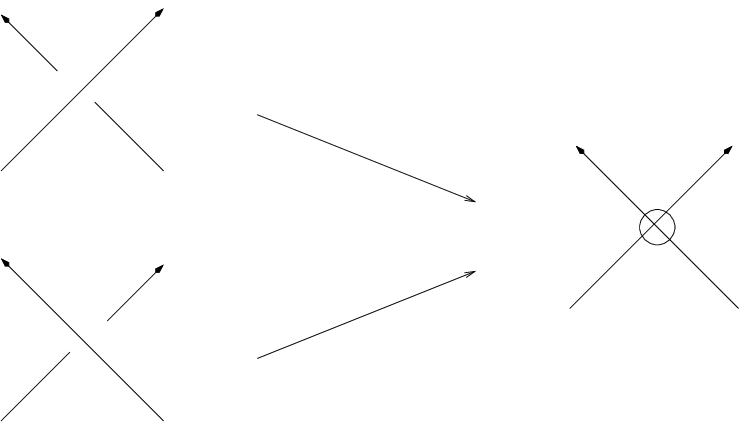}
\caption{\label{Fig:Resolutions2} The extended resolution.}
\end{figure}

  First, we digress briefly to extend our definition of $\overline{H}(\Gamma)$ to extended resolutions
$\Gamma'$. They are
resolutions, obtained by resolving a planar diagram of a link in one of the ways from
Figure~\ref{Fig:Resolutions} or as in Figure~\ref{Fig:Resolutions2}. To an extended resolution $\Gamma'$ we
associate a 2-complex $\overline{M}(\Gamma')$ in the same fashion as in subsection~\ref{Sbs:Hom}, we only have to add the description
of the factoriation, assiciated with the new resolution (the left part of Figure~\ref{Fig:Boundary2}). 
It is the tensor product of factorizations associated to arcs $x_3 x_1$ and $x_4 x_2$, that is the tensor
product of
$$\begin{CD} R @>{\pi_{13} - (n+1)\beta^n}>> R\lbrace 1-n\rbrace @>{x_3-x_1}>> R, \end{CD}$$
and
$$\begin{CD} R @>{\pi_{24} - (n+1)\beta^n}>> R\lbrace 1-n\rbrace @>{x_4-x_2}>> R, \end{CD}$$
where $\pi_{ij} = (x_i^{n+1}-x_j^{n+1})/(x_i-x_j)$.

We also define maps $\overline{\xi}_0:\overline{H}(\Gamma_1') \to \overline{H}(\Gamma_1), \overline{\xi}_1:
\overline{H}(\Gamma_1) \to \overline{H}(\Gamma_1')$ (Figure~\ref{Fig:Boundary2}). To do that, 
we use matrices $U_0, U_1, V_1, V_2$ from 
subsection~\ref{Sbs:hi1hi0}, with the roles of $x_3$ of $x_4$ exchanged. It therefore holds
\begin{equation}
\label{Eqn:xi0xi1}
 \overline{\xi}_0\overline{\xi}_1 = m(x_1) - m(x_3), \quad
 \overline{\xi}_1\overline{\xi}_0 = m(x_1) - m(x_3).
\end{equation}

We are now ready to produce a nonzero element ${\bf a}_\varphi\in Q_\varphi\overline{H}(\Gamma)$.
At each thick edge determine whether the state $\varphi$ is of type 1 or 2 (Figure~\ref{Fig:Patterns}). If it is
type 1, replace that thick edge with a $0$-resolution (the right to left trasformation of Figure~\ref{Fig:Boundary});
if it is type 2, replace that thick edge with the new resolution (the right to left transformation 
of Figure~\ref{Fig:Boundary2}). We thus obtain an extended resolution $\Gamma'$ which has only crossings
of the new type. The extended resolution $\Gamma'$ is a collection of circles (intersecting themselves at
new type crossings) $C_i, i=1,\ldots,k$. The state $\varphi$ induces a state $\varphi'$ of $\Gamma'$
in an obvious way; $\varphi'$ has the property that it evaluates to the same number at all edges of any
circle $C_i$, call that number $\varphi_i$. Note that $\varphi_i\not=\varphi_j$ for any two circles $C_i, C_j$
that either intersect or share, one arc each, the two arcs that resulted in performing the $0$-resolution
at $\Gamma$'s thick edge of type 1.

  Now, note the 2-complex $\overline{M}(\Gamma')$ is equal to the tensor product of 2-complexes (not just factorizations)
$\overline{M}_i(\Gamma')$, associated to circles $C_i$. The same must hold for the cohomology of $\overline{M}(\Gamma')$:
$\overline{H}(\Gamma')$ is the tensor product of cohomologies of $\overline{M}_i(\Gamma')$, each of them is
obviously equal to $\C[X_i]/(X_i^n=\beta^n)$. Therefore, we can choose a nonzero element in
$a_\varphi\in\overline{H}(\Gamma')$ that has $X_e a_\varphi=\beta\varphi'(e) a_\varphi$ for all thin edges $e$
of $\Gamma'$ - it is
$$ a_\varphi = \Pi_{i=1}^k (1+X_i/(\beta\varphi_i)+X_i^2/(\beta\varphi_i)^2+\ldots+X_i^{n-1}/(\beta\varphi_i)^{n-1})/n. $$
We define
$$ {\bf a}_\varphi := \Pi_{e_1} \overline{\chi}_0(e_1) \Pi_{e_2} \overline{\xi}_0(e_2) a_\varphi, $$
where the first product extends over thick edges $e_1$ of $\Gamma$ of type 1, and the second product
over thick edges $e_2$ of $\Gamma$ of type 2. First, since multiplication by edge variables commutes with
maps $\overline{\chi},\overline{\xi}$, we have $X_e {\bf a}_\varphi = \beta\varphi(e){\bf a}_\varphi$ and therefore
${\bf a}_\varphi=Q_\varphi{\bf a}_\varphi\in Q_\varphi\overline{H}(\Gamma)$.

Also, we claim ${\bf a}_\varphi$ is nonzero. This is because
$$ \Pi_{e_1} \overline{\chi}_1(e_1) \Pi_{e_2} \overline{\xi}_1(e_2) {\bf a}_\varphi = $$
$$ = \Pi_{e_1} \overline{\chi}_1(e_1)\overline{\chi}_0(e_1) \Pi_{e_2} \overline{\xi}_1(e_2)
\overline{\xi}_2(e_2) a_\varphi = $$
$$ = \Pi_{e_1} (\varphi_{i_{e_1}}-\varphi_{j_{e_1}}) 
     \Pi_{e_2} (\varphi_{i_{e_2}}-\varphi_{j_{e_2}}) a_\varphi $$
is a nonzero multiple of nonzero element $a_\varphi$. We used Equations (\ref{Eqn:hi0hi1}), (\ref{Eqn:xi0xi1}) and
the remark, made above, about a sufficient condition for $\varphi_i\not=\varphi_j$.

\end{proof}

  For each $\varphi\in S(\Gamma)$, choose a nonzero element of $Q_\varphi\overline{H}(\Gamma)$ 
 and call it ${\bf a}_{\varphi,\Gamma}$. According to the previous theorem, 
$\C {\bf a}_{\varphi,\Gamma} = Q_\varphi\overline{H}(\Gamma)$ and all vectors ${\bf a}_{\varphi,\Gamma}$ as
$\phi$ runs over all admissible states of $\Gamma$ form a basis for $\overline{H}(\Gamma)$.

  The boundary maps of $\overline{C}(L)$ are $\overline{\chi}_0,\overline{\chi}_1$. It is easy to determine
the action of these maps on basis vectors ${\bf a}$'s. Let $\Gamma_0, \Gamma_1$ be two resolutions, differing
locally as in Figure~\ref{Fig:Boundary}.
For any mark $i$ of $L$ we have
$$ m(x_i) \overline{\chi}_1 {\bf a}_{\varphi,\Gamma_1} = \overline{\chi}_1 m(x_i) {\bf a}_{\varphi,\Gamma_1} =
  \varphi(e_i) \overline{\chi}_1 {\bf a}_{\varphi,\Gamma_1}. $$
This implies that for any admissible state $\varphi\in S(\Gamma_1)$, the element $\overline{\chi}_1 {\bf a}_{\varphi,
\Gamma_1}$ is a multiple of ${\bf a}_{\varphi', \Gamma_0}$ where $\varphi'$ is the state
of $\Gamma_0$, obtained by retaining the same values of $\varphi$ at all marks. For states $\varphi$ of type 2,
the induced state $\varphi'$ is not well defined, since the two marks lying on one of the two arcs of $\Gamma_0$,
evaluate to different values. Therefore, for type 2 states $\varphi$ we have $\overline{\chi}_1 {\bf a}_{\varphi,
\Gamma_1}=0.$ For type 1 states $\varphi$, however, $\varphi'$ is an admissible state and we have
$\overline{\chi}_1 {\bf a}_{\varphi,\Gamma_1} = c {\bf a}_{\varphi',\Gamma_0}. $
Here, $c$ must be nonzero since 
$$ c\overline{\chi}_0{\bf a}_{\varphi',\Gamma_0} = \overline{\chi}_0\overline{\chi}_1{\bf a}_{\varphi,\Gamma_1} =
 (\varphi(e_{top,left}) - \varphi(e_{bottom,right})) {\bf a}_{\varphi,\Gamma_1} $$
is clearly nonzero.

  A similar calculation can be performed for $\overline{\chi}_0$. The behavior of boundary
maps on basis elements is summarized in Figure~\ref{Fig:Patterns}.

  We can now prove Theorem~\ref{Thm:Two}.
  
\begin{proof} (of Theorem~\ref{Thm:Two})

  As a $\C$-vector space, $\overline{C}(L)$ is spanned by elements ${\bf a}_{\varphi,\Gamma}$ for all resolutions
$\Gamma$ of $L$, and all admissible states of these resolutions. From the boundary map behavior 
(Figure~\ref{Fig:Patterns}) it is clear that the states that survive in the cohomology $\overline{H}(L)$
are precisely those, which are of type 2 at all crossings resolved to a 1-resolution, and of type 4 at all
crossings resolved to a 0-resolution. All such states clearly evaluate to the same value at every edge
in any component
of $L$, so are induced by a map
$$ \psi : \lbrace components\,\,of\,\,L\rbrace \to \Sigma_n. $$
Also, each such map $\psi$ defines precisely one admissible state: the one of a resolution, obtained by
resolving to $0$-resolution all crossings at which $\psi$-values of the two strands are equal, and
resolving to $1$-resolution all crossings at which $\psi$-values of the two strands differ.

To determine the cohomological degree of the state/resolution determined by some $\psi$, recall that
only $1$-resolutions contribute to the cohomological degree: the contribution is $1$ for positive
crossings and $-1$ for negative crossings.

\end{proof}

\bibliographystyle{plain}

\end{document}